\documentclass[10pt]{article}
\usepackage{amsmath,amssymb,amsthm}

\usepackage{enumitem}  

\usepackage{breqn}
\usepackage{setspace}

\usepackage{graphicx}	        
\usepackage{tikz-cd} 
\usepackage{bm}
\usepackage{xparse}
\usepackage{amsrefs}

\newtheorem*{question}{Question}

\allowdisplaybreaks[1]

\NewDocumentCommand{\evalat}{sO{\big}mm}{%
  \IfBooleanTF{#1}
   {\mleft. #3 \mright|_{#4}}
   {#3#2|_{#4}}%
}

\makeatletter
\newcommand{\colim@}[2]{%
  \vtop{\m@th\ialign{##\cr
    \hfil$#1\operator@font colim$\hfil\cr
    \noalign{\nointerlineskip\kern1.5\ex@}#2\cr
    \noalign{\nointerlineskip\kern-\ex@}\cr}}%
}
\newcommand{\colim}{%
  \mathop{\mathpalette\colim@{\rightarrowfill@\textstyle}}\nmlimits@
}
\makeatother

\makeatletter
\newcommand*{\rom}[1]{\expandafter\@slowromancap\romannumeral #1@}
\makeatother

\newcommand{\code}[1]{\texttt{#1}}

\usepackage[utf8]{inputenc}
\usepackage{listings}
\usepackage{dsfont}
\usepackage{float}
\usepackage{gensymb}
\usepackage{url}
\usepackage[toc,page]{appendix}

\begin{document}
\title{The Lagrangian Mechanics and Pseudo-Periodicity of The Many-Body Planar Pendulum System}

\author{Sergio Charles}

\date{10 March 2019}

\maketitle

\begin{abstract}
We present the Euler--Langrage equations for a many-body system of coupled planar pendulums. Hence, imposing initial condition data, the equations of motion are linearized and later developed in an idealized model for the pseudo-periodicity of the system as a function of the number of pendulums $N$. The result is empirically corroborated by comparing the model with data obtained via a numerical simulation, and by employing Kane's Method integrator in Python.
\end{abstract}

\tableofcontents

\newpage
\section{Introduction and Preliminaries\label{sec:intro}}
\fontsize{10.75}{12.1}\selectfont
A simple planar pendulum is a system with a bob of mass $m$ supported by a massless and inextensible rod of length $l$ which oscillates in two dimensions. Such a pendulum executes simple harmonic motion with regular oscillations. Up to small angle approximation, the perturbations are harmonic and can, thus, be described by sinusoidal functions. Historically, the approximate isochronism of the simple pendulum discovered by Galileo made it extremely useful in timekeeping \cite{nelson}. In fact, the simple pendulum was used by Newton as evidence for the universal law of gravitation. It provides a means of measuring gravitational acceleration $g\approx$9.8ms${}^{-2}$ near the surface of the Earth, which was used to postulate the spherical mass distribution of the Earth \cite{nelson}.

Lagrangian mechanics allows us to provide a comprehensive framework of physics and, thus, probe the underlying laws of nature. As an aspiring physicist, such an investigation is of great interest to me. While exploring problems in dynamics, I came across the interesting phenomenology of system resonance,  strange attractors, and topological mixing in dynamical systems with applications in various fields like signal processing and meteorology. Such research is of merit as the computational methods of Lagrangian mechanics manifest themselves in a wide range of fields, from the motion of planets in our solar system to geometrical optics and classical scalar field theory. The chaotic behavior exhibited by a system of pendulums is similarly seen in biological systems, self-assembly, the butter-fly effect, and engineering bridges.

\begin{figure}[H]
	\centering
	\includegraphics[width=1in]{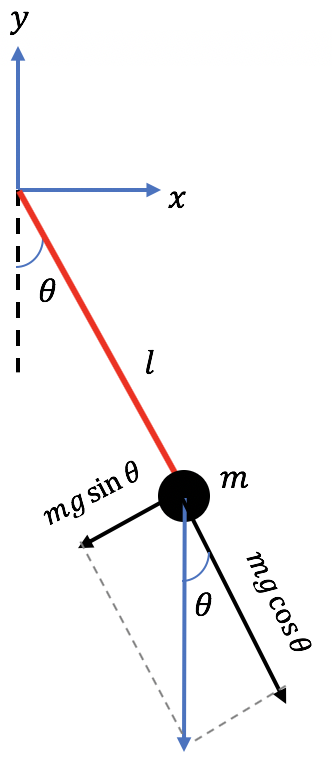}
	\caption
	{Schematic configuration of simple planar pendulum.}\label{simplependulum}
\end{figure}

Suppose the pendulum bob is displaced by an angle $\theta$ with respect to the vertical equilibrium of $\theta=0$ rad (see Figure \ref{simplependulum}). The mass will exert a restoring torque in the direction opposite to its displacement. In a first-order Taylor series expansion, the restoring torque is proportional to $\theta$, i.e. $\tau(\theta)=-\kappa\theta$ where $\kappa$ is the torsion constant. This constant is the restoring torque per unit angular displacement and the negative signature accounts for restoration \cite{shankar}. Thus, if $\theta$ is made to be positive in the counter-clockwise direction with respect to the positively-oriented vertical axis, then the torque will twist in the opposite direction. In fact, we find the torque to be $\tau=-mgl\sin\theta\approx-mgl\theta$ for sufficiently small $\theta$. 
As seen in the diagram, when we resolve the forces acting on the pendulum, the force due to gravity contributes to acceleration by $F=-mg\sin\theta=ma$ where $g\approx 9.8$ms${}^{-2}$ is the gravitational acceleration near the surface of the Earth. Thus, $F=-mg\sin\theta$ is the component of the gravitational force acting in the tangential direction opposite to the mass's motion. The arc length $s$ subtended by the mass bob, displaced by an angle $\theta$, is $s=l\theta$ such that the tangential velocity and acceleration are $v=\frac{ds}{dt}=l\frac{d\theta}{dt}$ and $a=\frac{d^2s}{dt^2}=l\frac{d^2\theta}{dt^2}$. 
Namely, the equation of motion for a pendulum attached to a mass $m$ by a rod reduces to:
\begin{equation}
\frac{d^2\theta}{dt^2}+\frac{g}{l}\sin\theta=0.
\end{equation} Suppose the initial amplitude of oscillation is $\theta(0):=\theta_0$. Assuming that the oscillation amplitude is sufficiently small, i.e. less than 1 radian, we use the small angle approximation $\sin\theta\approx\theta$ due to the Maclaurin series expansion of sine to find:
\begin{equation}\label{simpleode}
\frac{d^2\theta}{dt^2}+\frac{g}{l}\theta=0.
\end{equation} Thus, we make the ansatz that the solution is of the form $\theta(t)=A\cos(\omega t)$, with angular velocity $\omega$ and initial oscillation amplitude $A$, from which determine first and second-order time derivatives $\dot\theta(t)=-A\omega\sin(\omega t)$ and $\ddot\theta(t)=-A\omega^2\cos(\omega t)$. By substituting these derivatives into Equation \ref{simpleode}, we obtain:
\begin{equation}
A\cos(\omega t)\left(-\omega^2+\frac{g}{l}\right)=0
\end{equation} or that angular velocity is $\omega=\sqrt{g/l}$. Similarly, with initial conditions $(\theta(0),\dot\theta(0))=(\theta_0,0)$, we find that $\theta(t)=\theta_0\cos(\omega t)$ or
\begin{equation}
\theta(t)=\theta_0\cos\left(\sqrt{\frac{g}{l}}t\right).
\end{equation} Thus, the period describing the amount of time for a complete oscillation is 
\begin{equation}\label{simpleperiod}
T_0=\frac{2\pi}{\omega}=2\pi\sqrt{\frac{l}{g}}
\end{equation} for small $\theta_0$. The fact that the period is independent of the initial amplitude of oscillation $\theta_0$ means that the simple pendulum is isochronous. This property effectively means that the simple pendulum can be used to measure the acceleration due to gravity $g$. However, the error due to the approximation grows like $\theta^3$.

There are various physical assumptions of this idealized simple pendulum model, which we similarly adopt in the analysis of an $N$-body system \cite{shankar}:
\begin{enumerate}[label=(\roman*)]
\item Motion is constrained in two dimensions.
\item Friction due to air resistance is negligible.
\item The pendulum support is fixed.
\item The rod, from which the mass is attached, is massless and inextensible.
\end{enumerate}

For driven harmonic motion, pendulums will execute non-linear oscillations depending on the initial amplitude. We henceforth consider a system of $N$-pendulums attached by massless, inextensible rods of arbitrary mass and length. Later, we specialize the solution to the case for which the lengths and masses of the rods are equal.

\section{The Euler--Lagrange Equations for a System of $N$-Pendulums\label{sec:generalization}}
\begin{figure}[H]
	\centering
	\includegraphics[width=3.7in]{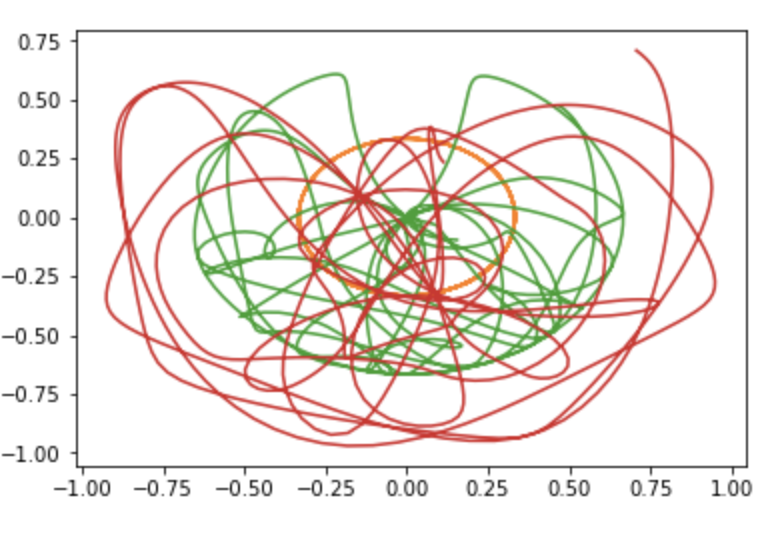}
	\caption
	{The chaotic trajectory of an $N=3$ dynamical system in phase space.}\label{attractor}
\end{figure}

The phase space trajectory of an $N$-body system of pendulums, where $N\ge 2$, is very distinct from that of a simple pendulum. That is, in the locality of small angle oscillation, the $N$-body pendulum consisting of $N$ masses exhibits beats, a characteristic of harmonic motion. As the total energy of the dynamical system increases, the oscillations behave chaotically. Insofar as an $N$-body problem can be described as a set of coupled ordinary differential equations, its chaotic response is an unexpected consequence of dynamical instability in Lorenz systems. That is, the real part of the eigenvalues of the system, when analytically solved, are less than unity, which therefore creates an attracting fixed point. As seen in Figure \ref{attractor}, the trajectory of the dynamical system in $(x,y)$ phase space exhibits chaotic behavior even for $N=3$ masses, which is obtained from the simulation, delineated in what follows. 

Thus, we seek to answer the following question and analyze the pseudo-periodicity of a simplified $N$-body system.
\begin{question}
What is the pseudo-time period of a system of $N$-coupled pendulums of unit $1$ m length and $1$ kg mass as a function of the number of pendulums $N$, all with initial oscillation amplitudes of $\pi/4$ radians and angular velocities of 0 rads${}^{-1}$?
\end{question}

We hence introduce the following independent, dependent, and controlled variables.
\begin{enumerate}[label=(\roman*)]
\item \textit{Independent variable:} The number of pendulums in the system, $N$.
\item \textit{Dependent variable:} The pseudo-oscillation time period of the dynamical system.
\item \textit{Controlled variable:} The length of the rods and the masses of the bobs, which are fixed to $1$ m and $1$ kg, respectively.
\end{enumerate}

\begin{figure}[H]
	\centering
	\includegraphics[width=1.6in]{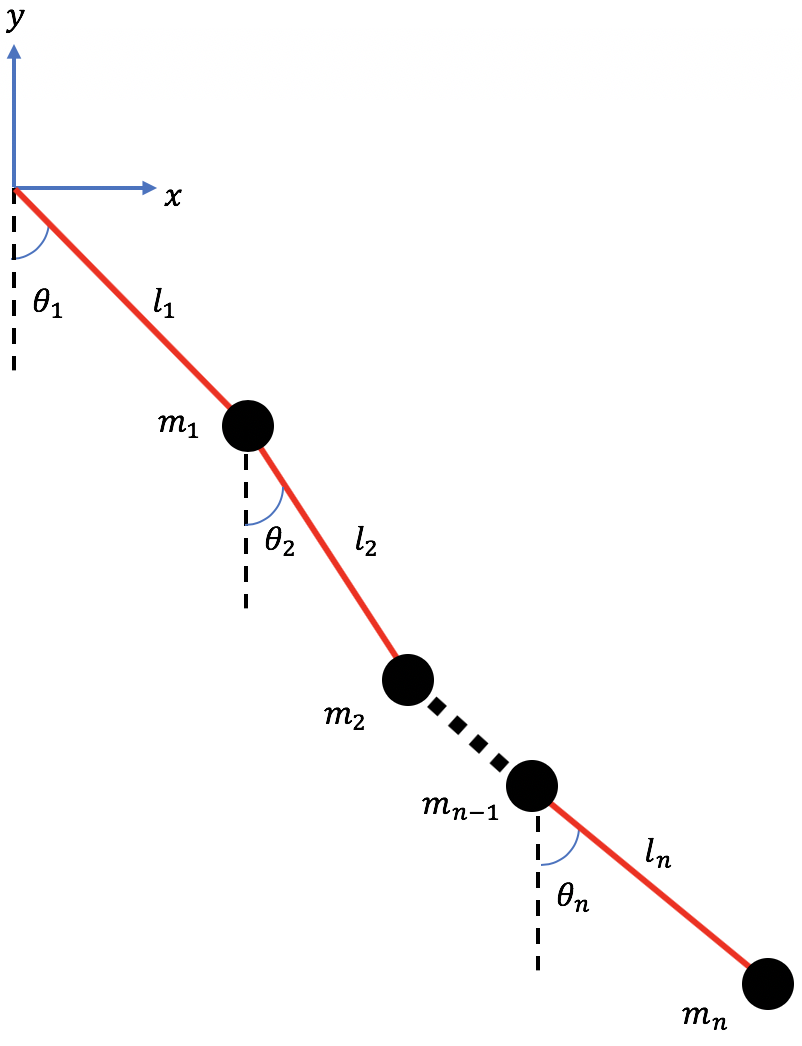}
	\caption
	{Free-body diagram--schematic configuration of $N$-body planar pendulum system.}\label{complexpendulum}
\end{figure}

We proceed by considering the Euler--Lagrange equations of motion for the dynamical system of coupled pendulums with arbitrary masses and lengths, which is completely characterized by generalized coordinates and velocities. In particular, the displacement angles from the vertical equilibrium of the pendulum $\theta_1,\dots,\theta_N$ and the corresponding angular velocities $\dot\theta_1,\dots,\dot\theta_N$ or $\omega_1,\dots,\omega_N$ are the set of generalized coordinates. Hence, we construct the Lagrangian for the double pendulum and numerically solve the system of Euler--Lagrange differential equations. The model of the $N$-body pendulum is illustrated in Figure \ref{complexpendulum}. 
The rods are assumed to be massless and frictionless with respective lengths of $l_1,\dots,l_N$. The point masses, represented by the attached balls of fixed radii, are $m_1,\dots,m_N$. We therefore impose a Euclidean coordinate system with the canonical metric in $\mathbb{R}^2$ where the origin $O$ is the point of suspension of the first pendulum of mass $m_1$. The coordinates are therefore given by $x_i=\sum_{k=1}^il_k\sin\theta_k$ and $y_i=-\sum_{k=1}^il_k\cos\theta_k.$ The kinetic energy $T(\dot\theta_1,\dots,\dot\theta_N)$ and the potential energy $V(\theta_1,\dots,\theta_N)$ are as follows:
\begin{equation}
\begin{split}
&T(\dot\theta_1,\dots,\dot\theta_N)=\frac{1}{2}\sum_{i=1}^N(\dot x_i^2+\dot y_i^2),\\
&V(\theta_1,\dots,\theta_N)=\sum_{i=1}^Nm_igy_i.
\end{split}
\end{equation} In particular, if $(x_i,y_i)$ denotes the Cartesian coordinate of the $i$-th mass $m_i$ of length $l_i$ of the pendulum for $1\le i\le N$, the Lagrangian formalism is
\begin{equation}
L(t,\theta_1,\dots,\theta_N,\dot\theta_1,\dots,\dot\theta_N)=T-V=\frac{1}{2}\sum_{i=1}^N(\dot x_i^2+\dot y_i^2)-\sum_{i=1}^Nm_igy_i
\end{equation} for $x_i=\sum_{k=1}^il_k\sin\theta_k$ and $y_i=-\sum_{k=1}^il_k\cos\theta_k.$ When $N=2$, this is simply 
\begin{equation}
\begin{split}
L=&\frac{1}{2}(m_1+m_2)l_1^2\dot\theta_1^2+\frac{m_2}{2}l_2^2\dot\theta_2^2+m_2l_1l_2\dot\theta_1\dot\theta_2\cos(\theta_1-\theta_2)\\&+(m_1+m_2)gl_1\cos\theta_1+m_2gl_2\cos\theta_2,
\end{split}
\end{equation} which, up to small angle approximation, is given by
\begin{equation}
\begin{split}
L=&\frac{1}{2}(m_1+m_2)l_1^2\dot\theta_1^2+\frac{m_2}{2}l_2^2\dot\theta_2^2+m_2l_1l_2\dot\theta_1\dot\theta_2\\&-\frac{1}{2}(m_1+m_2)gl_1\theta_1^2+\frac{m_2}{2}gl_2\theta_2^2.
\end{split}
\end{equation} As such, we determine the Euler--Lagrange equation for each angle via
\begin{equation}
\frac{\partial L}{\partial\theta_i}=\frac{d}{dt}\left(\frac{\partial L}{\partial \dot\theta_i}\right)
\end{equation} for $1\le i\le N$, which yields a system of second-order non-linear differential equations. 
We compute the first and second-order time derivatives of the position $(x_i,y_i)$. In particular,
\begin{equation}
\begin{split}
&\dot x_i=\sum_{k=1}^il_k\dot\theta_k\cos\theta_k,\hspace{0.9 cm} \ddot x_i=\sum_{k=1}^il_k(\ddot\theta_k\cos\theta_k-\dot\theta_k^2\sin\theta_k)\\
&\dot y_i=\sum_{k=1}^il_k\dot\theta_k\sin\theta_k,\hspace{1 cm} \ddot y_i=\sum_{k=1}^il_k(\ddot\theta_k\sin\theta_k+\dot\theta_k^2\cos\theta_k).
\end{split}
\end{equation}
Applying the multivariable chain-rule, we find:
\begin{equation}\label{leftlagrangian}
\begin{split}
\frac{\partial L}{\partial\theta_i}&=\frac{\partial}{\partial\theta_i}\sum_{j=1}^N\frac{1}{2}m_j(\dot x_j^2+\dot y_j^2)-gm_jy_j\\&
=\sum_{j=1}^N\frac{1}{2}m_j\left[2\dot x_j\frac{\partial \dot x_j}{\partial\theta_i}+2y_j\frac{\partial \dot y_j}{\partial\theta_i}\right]-gm_j\frac{\partial y_j}{\partial\theta_i}\\&=\sum_{j=1}^Nm_j\left[\dot x_j\frac{\partial \dot x_j}{\partial\theta_i}+\dot y_j\frac{\partial \dot y_j}{\partial\theta_i}-g\frac{\partial y_j}{\partial\theta_i}\right].
\end{split}
\end{equation} Similarly,
\begin{equation}\label{rightlagrangian}
\begin{split}
\frac{\partial L}{\partial\dot\theta_i}=&\frac{1}{2}\sum_{j=1}^Nm_j\left[2\dot x_j\frac{\partial \dot x_j}{\partial\dot\theta_i}+2\dot y_j\frac{\partial\dot y_j}{\partial\dot\theta_i}\right]-g\sum_{i=1}^Nm_j\frac{\partial y_j}{\partial\dot\theta_i}\\
&=\sum_{i=1}^Nm_j\left[\dot x_j\frac{\partial\dot x_j}{\partial\dot\theta_i}+\dot y_j\frac{\partial\dot y_j}{\partial\dot\theta_i}-g\frac{\partial\dot y_j}{\partial\dot\theta_i}\right].
\end{split}
\end{equation}
Therefore, the Euler--Lagrange equations 
\begin{equation*}
\frac{\partial L}{\partial\theta_i}=\frac{d}{dt}\left(\frac{\partial L}{\partial\dot\theta_i}\right)
\end{equation*} become (see Appendix \ref{sec:appendixa} for calculation):
\begin{equation}
\begin{split}
&-g\sin\theta_i\sum_{j=i}^Nm_j-\dot\theta_i\sum_{j=1}^Nl_j\dot\theta_j\sin(\theta_i-\theta_j)\sum_{j=1}^Nl_j\dot\theta_j\sin(\theta_i-\theta_j)\sum_{k=\pi_{ij}}^Nm_k\\
&=\sum_{j=1}^N
l_j\left[\ddot\theta_j\cos(\theta_i-\theta_j)+(\dot\theta_j^2-\dot\theta_i\dot\theta_j)\sin(\theta_i-\theta_j)\right]\sum_{k=\pi_{ij}}^Nm_k.
\end{split}
\end{equation}
This is equivalent to the standard LU decomposition matrix equation:
\begin{equation}\sum_{j=1}^N\textbf{A}_{ij}\ddot\theta_j=\textbf{b}_i
\end{equation} with
\begin{equation}
\begin{split}
\textbf{A}_{ij}&=l_j\cos(\theta_i-\theta_j)\sum_{k=\pi_{ij}}^Nm_k,\\
\textbf{b}_i&=-g\sin\theta_i\sum_{j=i}^Nm_j-\sum_{j=1}^Nl_j\left[\dot\theta_j^2\sin(\theta_i-\theta_j)\right]\sum_{k=\pi_{ij}}^Nm_k,\\
\pi_{ij}&=\begin{cases}i, &\text{if }i\le j\\
j, & \text{else.}\end{cases}
\end{split}
\end{equation} 

\section{Linearization and Pseudo-Periodicity\label{sec:linearization}}
The equations of motion for the system of $N$-pendulums are:
\begin{equation}\label{eom}
\sum_{j=1}^N\left[l_j\cos(\theta_i-\theta_j)\sum_{k=\pi_{ij}}^Nm_k\ddot{\theta}_j\right]=-g\sin\theta_i\sum_{j=i}^Nm_j-\sum_{j=1}^Nl_j\left[\dot\theta_j^2\sin(\theta_i-\theta_j)\right]\sum_{k=\pi_{ij}}^Nm_k.
\end{equation} This is a system of non-linear ordinary differential equations with no analytic solution. Up to small angle approximation, the oscillations of the pendulums near the gravity-induced zero equilibrium $(\theta_1,\dots,\theta_n)^T=\mathbf{0}$ are uniquely determined by a linearized set of ordinary differential equations. Let 
\begin{equation} 
\bf{\Theta}:=\begin{bmatrix}\theta_1\\ \vdots \\ \theta_n\end{bmatrix},
\end{equation}
\begin{equation}
\bf M=\begin{bmatrix}l_1\sum_{k=\pi_{11}}^Nm_k & & & \bf 0\\
 & l_2\sum_{k=\pi_{22}}^Nm_k & & & \\ & & \ddots & \\
\bf 0 & & & l_n\sum_{k=\pi_{NN}}^Nm_k\end{bmatrix},
\end{equation}
\begin{equation}
\bf K=\begin{bmatrix}l_1\sum_{k=\pi_{11}}^Nm_k & & \bf 0\\ & \ddots & \\ \bf 0 & & l_n\sum_{k=\pi_{NN}}^Nm_k\end{bmatrix},
\end{equation} and
\begin{equation}
\bf L=\begin{bmatrix}gm_1 & & \bf 0\\ & \ddots & \\ \bf 0 & & gm_N\end{bmatrix}.
\end{equation} Thus, the equations of motion are $\textbf{M}\ddot{ \mathbf{\Theta}}+\mathbf{K\dot\Theta}^T\mathbf{\dot\Theta}+\mathbf{L\Theta}=\mathbf{0}$ (see Appendix \ref{sec:appendixlin}). For a simple pendulum, this equation describes a freely-evolving undamped system with a natural frequency \cite{math24}. In the case of the $N$-body system, the solution will consist of oscillations of $N$ characteristic frequencies, i.e. normal modes. The normal modes are the real part of the complex-valued vector function 
\begin{equation} \mathbf{\Theta}(t)=\begin{bmatrix} \theta_1(t)\\\vdots\\\theta_n(t)\end{bmatrix}=\text{Re}\left(\begin{bmatrix}\mathbf{H}_1\\\vdots\\ \mathbf{H}_n\end{bmatrix}e^{i\omega t}\right)
\end{equation} where $\mathbf{H}_1,\dots\mathbf{H}_n$ are eigenvectors and $\omega$ is the real frequency of the entire system \cite{math24}. To first-order approximation, we have ${\dot{\mathbf{\Theta}}}^T{\dot{\mathbf{\Theta}}}\approx \mathbf{0}$ so the equations of motion become $\textbf{M}\ddot{\mathbf{\Theta}}+\mathbf{L\Theta}=\mathbf{0}$. This system can similarly be linearized using a Jacobian. The normal frequencies are thereby determined as solutions to the auxiliary determinant equation:
\begin{equation}
\det(\textbf{K}-\omega^2\textbf{M})=0.
\end{equation} That is, if $\omega_j$ is the angular velocity for $m_j$,
\begin{equation}
\det(\textbf{K}-\omega^2\textbf{M})=\prod_{j=1}^N\left[gm_j-\omega_j^2l_j\sum_{k=\pi_{ij}}^Nm_k\right]=0
\end{equation} because $\mathbf{K}-\omega^2\mathbf{M}$ is a diagonal matrix. Therefore, $\omega_j\approx \sqrt{\frac{gm_j}{l_j M}}$ where $M=\sum_{i=1}^Nm_i\approx \sum_{k=\pi_{ij}}^Nm_k$ so the average angular velocity of the system response is $\overline{\omega}(N)=\frac{1}{N}\sum_{j=1}^N\sqrt{\frac{gm_j}{l_jM}}$ and the pseudo-time period as a function of the number of pendulums $N$ is approximately 
\begin{equation}\label{npesudoperiod}
\begin{split}
T(N)&\approx\frac{2\pi N}{\sqrt{\frac{gm_1}{l_1M}}+\dots+\sqrt{\frac{gm_N}{l_NM}}}\\&=\frac{2\pi N}{\sum_{i=1}^N\sqrt{\frac{gm_i}{l_iM}}}.
\end{split}
\end{equation} If all pendulums have rods of length $l$ and bobs of equal mass $m$, the pseudo-period becomes 
\begin{equation}
T(N)\approx\frac{2\pi N}{N\sqrt{\frac{gm}{lmn}}}=2\pi\sqrt{\frac{l}{g}}\sqrt{N}.
\end{equation} Comparing the pseudo-period for a system of $N$-pendulums of equal masses $m$ and lengths $l$ with that of a simple pendulum in Equation \ref{simpleperiod}, it is evident that the period scales by $O(\sqrt{N}).$

\section{Higher-Order Term Corrections to Pseudo-Period and Fractional Uncertainty}
We remark that the linearized equations of motion in Equation \ref{eom} are second-order ordinary linear differential equations. For an initial value problem, by the existence and uniqueness Picard--Lindel\"{o}f theorem for differential equations, the solution requires $2N$ initial conditions for initial angular displacement and initial angular velocity. In particular, we must specify $\mathbf{\Theta}(0):=\mathbf{\Theta}_0$ and $\dot{\mathbf{\Theta}}(0):=\dot{\mathbf{\Theta}_0}$. Accounting for higher-order corrections, the exact solution to Equation \ref{eom} involves a Jacobian elliptic sine function \cite{nelson}. Thus, the period can be written in terms of an elliptic integral of first kind in a perturbative expansion. For a real system of $N$ pendulums, the period of oscillation $T$ is $T_0+\Delta T=T_0(1+\Delta T/T_0)$ where $T_0$ is the idealized period from Equation \ref{npesudoperiod} and $\Delta T/T_0$ is fractional random uncertainty in the experiment. If $\theta_0$ is the initial angular displacement, the higher-order term correction to the period is \cite{nelson}:
\begin{equation}
\begin{split}\label{absolute}
\frac{\Delta T}{T_0}&=\sum_{n=1}^{\infty}\left(\frac{(2n)!}{2^{2\pi}(n!)^2}\right)^2\sin^{2n}\left(\frac{\theta_0}{2}\right)\\
&=\frac{1}{16}\theta_0^2+\frac{11}{3072}\theta_0^4+\cdots.
\end{split}
\end{equation} When $\theta_0\approx \pi/4$ rad for the initial value problem, we find
\begin{equation}\label{45error}
\frac{\Delta T}{T_0}=\sum_{n=1}^{\infty}\left(\frac{(2n)!}{2^{2\pi}(n!)^2}\right)^2\sin^{2n}\left(\frac{\pi/4}{2}\right)\approx 0.049.
\end{equation}
Therefore, the exact pseudo-period solution for a general system is given by
\begin{equation}\label{periodexact}
T(N,\theta_0)=\frac{2\pi N}{\sum_{i=1}^N\sqrt{\frac{gm_i}{l_iM}}}\Bigg[1+\sum_{n=1}^{\infty}\left(\frac{(2n)!}{2^{2\pi}(n!)^2}\right)^2\sin^{2n}\left(\frac{\theta_0}{2}\right)\Bigg].
\end{equation} The difference between the real period and ideal period for $\theta\ll 1$ is said to be the \textit{circular error}.

\section{Procedure and Analysis of Kane's Method Simulation for Chaotic Dynamical System of Pendulums\label{sec:approximation}}

Kane's method provides a first-order numerical scheme used to approximate non-linear differential equations of motion. We can use Kane's method by importing the \code{mechanics} library from the \code{Sympy} package in Python to integrate the equations of motion for a system of $N$ pendulums with arbitrary masses and lengths \cite{sympy}. However, as we will see, such a method will introduce local and global truncation errors in the numerical approximation.

\begin{figure}[H]
	\centering
	\includegraphics[width=3.7in]{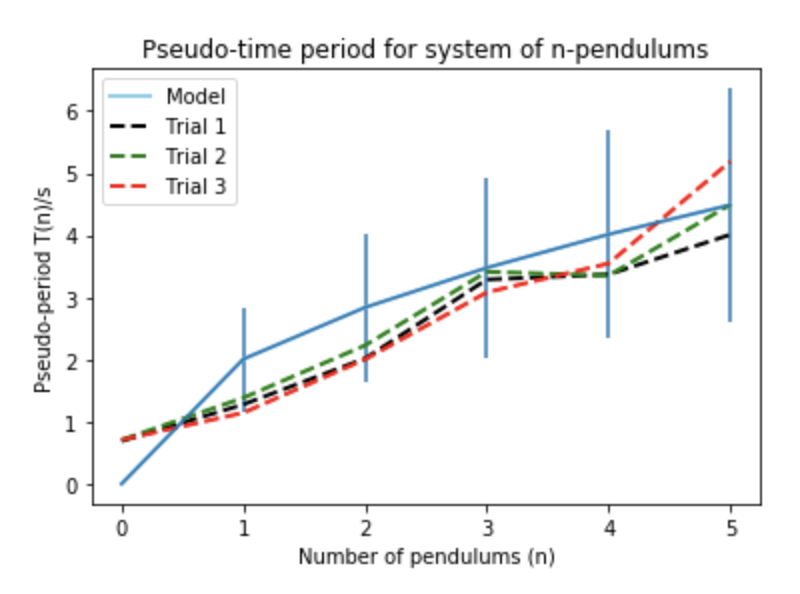}
	\caption
	{Simulated trials of pseudo-time period for $N=5$ with absolute error bars.}\label{fig1}
\end{figure}

\begin{figure}[H]
	\centering
	\includegraphics[width=3.56in]{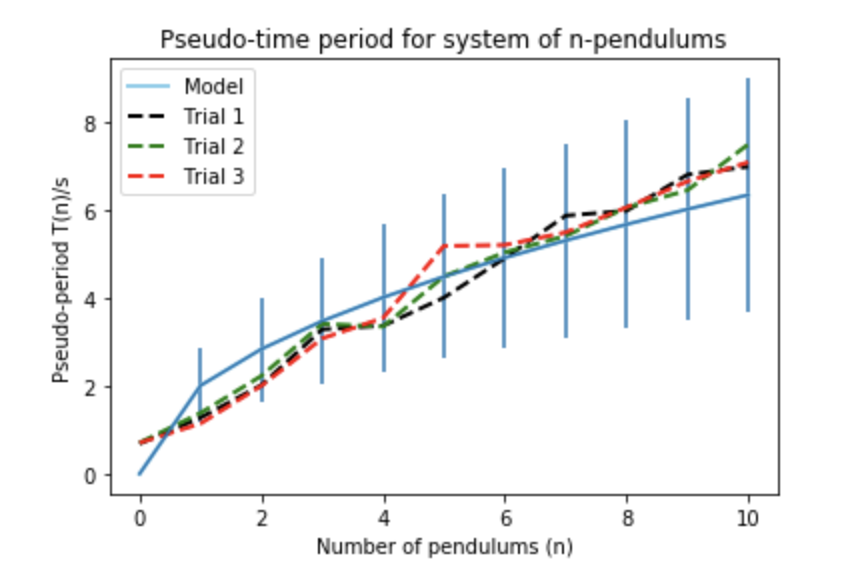}
	\caption
	{Simulated trials of pseudo-time period for $N=10$ with absolute error bars.}\label{fig2}
\end{figure}

\begin{figure}[H]
	\centering
	\includegraphics[width=3.51in]{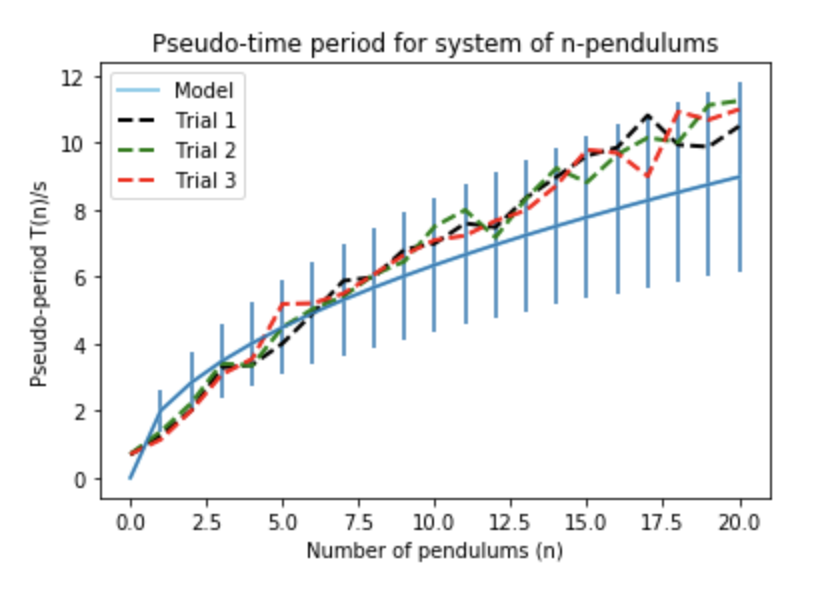}
	\caption
	{Simulated trials of pseudo-time period for $N=20$ with absolute error bars.}\label{fig3}
\end{figure}

\begin{figure}[H]
	\centering
	\includegraphics[width=3.77in]{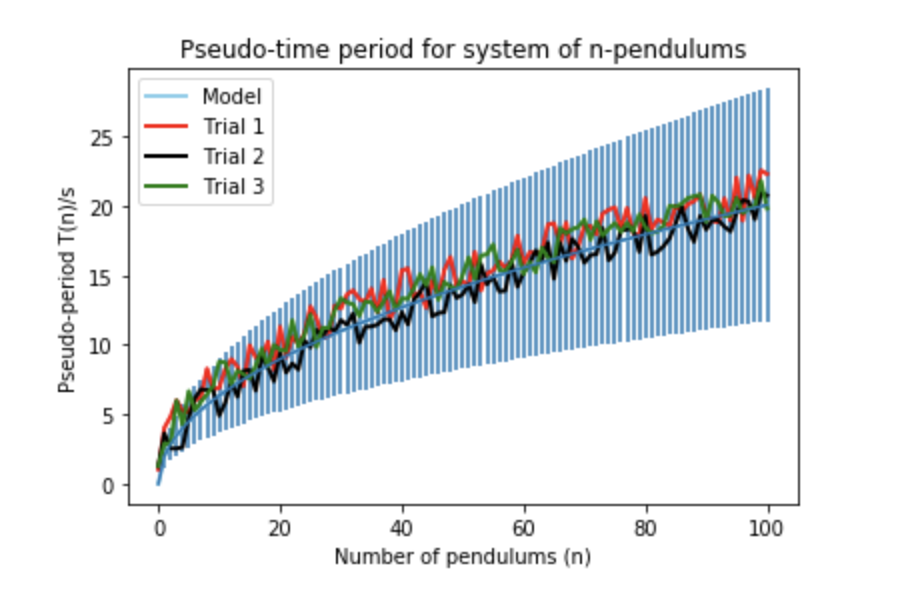}
	\caption
	{Simulated trials of pseudo-time period for $N=100$ with absolute error bars.}\label{fig4}
\end{figure}

The initial value problem under consideration is a system with initial oscillation amplitudes of $(\theta_1(0),\dots,\theta_N(0))=(\pi/4,\dots,\pi/4)$ and initial angular velocities of $(\dot\theta_1(0),\dots,\theta_N'(0))=(0,\dots,0).$ The procedure is as follows:
\begin{enumerate}[label=(\roman*)]
\item Repeat steps (ii) to (vi) for $N=5$, $10$, $20,$ and $100$.
\item Displace all $N$ bobs of the system by an angle of $\theta=\pi/4$ radians from the vertical equilibrium position of $\theta=0$, bearing in mind zero offset error. 
\item Release all bobs with zero velocity from the initial angular position. Exercise caution when handling large bobs if conducting the experiment for masses other than $1$ kg.
\item Using a stopwatch, measure the amount of time in seconds each bob takes to complete an entire oscillation, whereby a period is said to be \textit{pseudo-complete} if the angular velocity of the bob changes from zero to negative, to positive, and back to zero or the bob returns to its original vertical position (not necessarily horizontal). That is, the bob must change its direction twice or return to its original vertical height. This definition only holds true for small $\theta$.
\item Determine the average of the times taken by each bob and report this as the measured pseudo-period.
\item Iterate over steps (ii) to (v) for $t:=3$ trials to minimize systematic error $e$ and standard deviation $\sigma=e\sqrt{t}$. 
\end{enumerate}
When we run the simulation, in order to obtain different trials, we perturb the initial angular displacements by an amount within the range $[0,0.017]$ rad of absolute uncertainty in radian measurement. Thus, the probability of having trials with the same initial oscillation amplitude is effectively zero. Namely, for the initial angular displacement of $\mathbf{\Theta}_0=\mathds{1}^T\theta_0$, we perturb initial condition $\theta_0=\pi/4$ by $\delta\theta_0$ to obtain the new trial equation of
\begin{equation}
T(N,\theta_0+\delta\theta_0)=\frac{2\pi N}{\sum_{i=1}^N\sqrt{\frac{gm_i}{l_iM}}}\Bigg[1+\sum_{n=1}^{\infty}\left(\frac{(2n)!}{2^{2\pi}(n!)^2}\right)^2\sin^{2n}\left(\frac{\theta_0+\delta\theta_0}{2}\right)\Bigg].
\end{equation} We note that the simulation will run for $10$ seconds over $1000$ frames and, thus, find average periods for each mass over said time. See Appendix \ref{sec:appendix} for the Python implementation of the \code{KanesMethod(...)} to integrate the Euler--Lagrange equations of motion, which is called via the \code{Simulation()} method. We run the simulation for $N=5$, $10$, $20$, and $100$ to obtain the following graphs shown in Figure \ref{fig1}, \ref{fig2}, \ref{fig3}, and \ref{fig4}.

Thus, according to each simulated trial, the pseudo-periodicity varies stochastically; however, it does so with respect to the central mean of the model described in Section \ref{sec:generalization}. It is evident that the errors introduced in the linear approximations used to analytically solve the Euler--Lagrange equations have propagated throughout the calculation. That is, the determination of the pseudo-period depends on the nature of determinant or characteristic equation obtained after linearizing the set of non-linear Euler--Lagrange equations. It is empirically clear that the degree of accuracy of Kane's method for this numerical approximation scheme decreases as the number of pendulums increases, i.e., the global truncation error increases as $N\to\infty$. 

Given that all data, with the exception of one point, lies within the threshold of uncertainty, we can reasonably infer that the model developed is suitable for small angles. We now determine the (cumulative) average decimal percent error $e$ of the simulation from the model for $N=5$, $10$, $20$, and $100$ and each trial. In particular, for $T_i$ the model prediction and for $S_i$ the simulated results for each $1\le i\le N$, we find $\frac{|S_i-T_i|}{S_i}=\frac{|\text{real period}-\text{ideal period}|}{\text{real period}}$ and then compute the approximate error 
\begin{equation*}
e\approx\frac{1}{N}\sum_{i=1}^N\frac{|S_i-T_i|}{S_i}.
\end{equation*}

\begin{table}[H]
\centering
\resizebox{13cm}{.8cm}{
    \begin{tabular}{ | c | c | c | c | c | }
    \hline
    Trial 1 & $N=5$ & $N=10$ & $N=20$ & $N=100$ \\ \hline
    Initial oscillation amplitude [rad] & 0.77578 rad & 0.797476 rad & 0.78477 rad & 0.777 rad \\ \hline
    Cumulative average decimal error & 0.142 & 0.108 & 0.147 & 0.376  \\ \hline
    \end{tabular}
    }
\caption {Trial 1 average decimal errors.}\label{table1}
\end{table} 

\begin{table}[H]
\centering
\resizebox{13cm}{.8cm}{
    \begin{tabular}{ | c | c | c | c | c | }
    \hline
    Trial 2 & $N=5$ & $N=10$ & $N=20$ & $N=100$ \\ \hline
    Initial oscillation amplitude [rad] & 0.77253 rad & 0.77720 rad & 0.79765 rad & 0.77676 rad \\ \hline
    Cumulative average decimal error & 0.194 & 0.137 & 0.158 & 0.207  \\ \hline
    \end{tabular}
    }
\caption {Trial 2 average decimal errors.}\label{table2}
\end{table} 

\begin{table}[H]
\centering
\resizebox{13cm}{.8cm}{
    \begin{tabular}{ | c | c | c | c | c | }
    \hline
    Trial 3 & $N=5$ & $N=10$ & $N=20$ & $N=100$ \\ \hline
    Initial oscillation amplitude [rad] & 0.77419 rad & 0.77749 rad & 0.78107 rad & 0.78599 rad \\ \hline
    Cumulative average decimal error & 0.222 & 0.149 & 0.161 & 0.382  \\ \hline
    \end{tabular}
    }
\caption {Trial 3 average decimal errors.}\label{table3}
\end{table} 
 
According to Tables \ref{table1}, \ref{table2}, \ref{table3}, it is apparent that the local truncation error is relatively well-behaved for a small number of pendulums, i.e. the decimal error is under $0.23$ for all $N=5$ observations. However, for $N=100$, the decimal error grows exponentially, as is the case in Trial 1 with average error $0.142$ for $N=5$ and an average error of $0.376$ for $N=100$. This large increase in local truncation error may be attributed to the chaotic behavior of the large system of pendulums. In particular, as more pendulums are added, the more the system exhibits beats and thus behaves stochastically, which forces the error to rise artificially high.

When we implement Kane's method, we must first transform the equations of motion
\begin{equation}
\ddot{\mathbf{\Theta}}(t)=\mathbf{F}(t,\mathbf{\Theta}(t),\dot{\mathbf{\Theta}}(t))=-\mathbf{K}{\dot{\mathbf{\Theta}}}^T\dot{\mathbf{\Theta}}-\mathbf{L}\mathbf{\Theta}
\end{equation} into a first-order system of ODEs. We introduce the auxiliary matrices $\bm{\alpha}_1(t)=\mathbf{\Theta}(t)$ and $\bm{\alpha}_2(t)=\dot{\mathbf{\Theta}}(t).$ Thus, we obtain a first-order system in $\bm{\alpha}(t)$:
\begin{equation}\label{system}
\dot{\bm{\alpha}}(t)=\begin{bmatrix}
           \dot{\bm{\alpha}}_1(t)\\
           \dot{\bm{\alpha}}_2(t)
         \end{bmatrix}
         =\begin{bmatrix}
           \dot{\mathbf{\Theta}}(t)\\
           \ddot{\mathbf{\Theta}}(t)
         \end{bmatrix}
         =\begin{bmatrix}
           \bm{\alpha}_2(t)\\
           \mathbf{F}(t,\mathbf{\Theta},\dot{\mathbf{\Theta}})
         \end{bmatrix}.
\end{equation}

We can use the 4th order Runge-Kutta method (RK4) to numerically solve this system, which is similar to the 2nd order Kane's method (see Appendix \ref{sec:appendixrk4}). More precisely, RK4 is a fourth-order method and, thus, the local truncation error, defined as error per step $\tau$, is proportional to the fifth power of the step $O(\tau^5)$. On the other hand, global truncation error, defined to be the error at a given time, has order $O(\tau^4)$ \cite{math24}. Thus, if we halve the step size, the global error is reduced by a factor of $1/16$.

Consider the Taylor series expansion of the matrix $\mathbf{\Theta}(t)$ around $t_0$: 
\begin{equation}
\mathbf{\Theta}(t_0+\tau)=\mathbf{\Theta}(t_0)+\tau\dot{\mathbf{\Theta}}(t_0)+\frac{1}{2}\tau^2\ddot{\mathbf{\Theta}}(t_0)+O(\tau^3).
\end{equation} The local truncation error manifests itself after a single step and is said to be the difference between the numerical solution $\mathbf{\Theta}_1$ after one step $\tau$ and the the exact solution to Equation \ref{eom} at $t_1=t_0+\tau$. Thus, the local truncation error $e$ is the difference of $\mathbf{\Theta}_1=\mathbf{\Theta}_0+\tau\mathbf{F}(t_0,\mathbf{\Theta}_0)$ and the aforementioned Taylor series approximation:
\begin{equation}
e=\|\mathbf{\Theta}_1-\mathbf{\Theta}(t_0+\tau)\|=-\frac{1}{2}\tau^2\|\ddot{\mathbf{\Theta}}(t_0)\|+O(\tau^3).
\end{equation} In general, at $t=t_n$, the local truncation error is:
\begin{equation}
e=\|\mathbf{\Theta}_n-\mathbf{\Theta}(t_{n-1}+\tau)\|=-\frac{1}{2}\tau^2\|\ddot{\mathbf{\Theta}}(t_{n-1})\|+O(\tau^3).
\end{equation}
As such, the \textit{sign} of the truncation error depends on the sign of local curvature, i.e. concavity and convexity due to the second derivative, of the integral curve satisfying Equation \ref{eom}. Thus, the real degree of truncation error in the amplitude obtained in Tables \ref{table1}, \ref{table2}, and \ref{table3} can be analytically bounded by determining the circular error between the real and ideal amplitude. That is, since $T=T_0\pm\Delta T=T_0\pm 0.049T_0$, let
\begin{equation}
\begin{split}
\hat{e}:&=1-\frac{|\text{real period}-\text{ideal period}|}{\text{real period}}\\
&=1-\frac{|T(N,\theta_0)-T_0(N)|}{T(N,\theta_0)}\\
&=\frac{1}{1+\sum_{n=1}^{\infty}\left(\frac{(2n)!}{2^{2\pi}(n!)^2}\right)^2\sin^{2n}\left(\frac{\theta_0}{2}\right)},
\end{split}
\end{equation} which is clearly independent of $N$ so we do not need to take an average. For $\theta_0\approx \pi/4$ rad, we find
\begin{equation}
\hat{e}=\frac{1}{1+\sum_{n=1}^{\infty}\left(\frac{(2n)!}{2^{2\pi}(n!)^2}\right)^2\sin^{2n}\left(\frac{\pi/4}{2}\right)}\approx 0.309
\end{equation} up to order $20$ in the expansion. In fact, when we evaluate the entire series, it diverges so the error $\hat e$ effectively tends to $0$. Therefore, since the majority of empirically observed average circular errors for each trial in Tables \ref{table1}, \ref{table2}, and \ref{table3}, with a few outliers, lie within this analytic threshold, the model is deemed accurate. We compute the real period $T=T_0\pm\Delta T$ where $T_0=\sqrt{lN/g}=\sqrt{N/g}$. To find an upper-bound, rather than taking an average, we calculate the ideal period at the maximum value for each trial of $N=5$, $10$, $20$, and $100$. Using a sixth order expansion of Equation \ref{absolute} and $\theta_0\approx \pi/4$ rad, the approximate fractional uncertainty of the period is $\Delta T/T_0\approx 0.311$, from which we obtain:

\begin{table}[H]
\centering
\resizebox{\textwidth}{!}{
    \begin{tabular}{ | c | c | c | c | c | }
    \hline
    Trial 1 & Ideal period $T_0=2\pi\sqrt{N/g}$ [s] & Absolute uncertainty $\Delta T=0.311T_0$ [s] & Real period $T=T_0\pm\Delta T$ [s] \\ \hline
    $N=5$ & $4.487$ & $1.385$ & $4.487\pm 1.385$ s  \\ \hline
    $N=10$ & $6.347$ & $1.973$ & $6.347\pm1.973$  \\ \hline
    $N=20$ & $8.976$ & $2.792$ & $8.976\pm2.792$ \\ \hline
    $N=100$ & $20.071$ & $6.242$ & $20.071\pm6.242$ \\ \hline
    \end{tabular}
    }
\caption {Trial 1--Real and ideal period with absolute uncertainty.}\label{table_1}
\end{table} 

\begin{table}[H]
\centering
\resizebox{\textwidth}{!}{
    \begin{tabular}{ | c | c | c | c | c | }
    \hline
    Trial 2 & Ideal period $T_0=2\pi\sqrt{N/g}$ [s] & Absolute uncertainty $\Delta T=0.427T_0$ [s] & Real period $T=T_0\pm\Delta T$ [s] \\ \hline
    $N=5$ & $4.487$ & $1.916$ & $4.487\pm 1.916$ s  \\ \hline
    $N=10$ & $6.347$ & $2.710$ & $6.347\pm2.710$  \\ \hline
    $N=20$ & $8.976$ & $3.832$ & $8.976\pm3.832$ \\ \hline
    $N=100$ & $20.071$ & $8.570$ & $20.071\pm8.570$ \\ \hline
    \end{tabular}
    }
\caption {Trial 2--Real and ideal period with absolute uncertainty.}\label{table_2}
\end{table}

\begin{table}[H]
\centering
\resizebox{\textwidth}{!}{
    \begin{tabular}{ | c | c | c | c | c | }
    \hline
    Trial 3 & Ideal period $T_0=2\pi\sqrt{N/g}$ [s] & Absolute uncertainty $\Delta T=0.293T_0$ [s] & Real period $T=T_0\pm\Delta T$ [s] \\ \hline
    $N=5$ & $4.487$ & $1.315$ & $4.487\pm 1.315$ s  \\ \hline
    $N=10$ & $6.347$ & $1.860$ & $6.347\pm1.860$  \\ \hline
    $N=20$ & $8.976$ & $2.630$ & $8.976\pm2.630$ \\ \hline
    $N=100$ & $20.071$ & $5.880$ & $20.071\pm5.880$ \\ \hline
    \end{tabular}
    }
\caption {Trial 3--Real and ideal period with absolute uncertainty.}\label{table_3}
\end{table}

We remark that the numerical approximation will under-estimate or over-estimate the system response for small $N$ because we are approximating $M\approx\sum_{k=\pi_{ij}}^Nm_k$. Thus, the solution will more closely approximate the exact response as $N\to\infty$.

\section{Conclusion and Improvements\label{sec:conclusion}}
We have mathematically derived and empirically verified the functional relationship between the pseudo-periodicity of a system of $N$-coupled pendulums and the number of attached pendulums, $N$. In particular, by linearizing the Euler--Lagrange equations, we were able to analytically develop a model for the ideal period $T_0$, which is accurate up to third order in a Taylor series. Moreover, introducing correction terms $\Delta T/T_0$, we explicitly determined the pseudo-period for an initial value problem where all periods have mass $1$ kg, length $1$ m, initial oscillation amplitude $\pi/4$ rad, and angular velocity $0$rads${}^{-1}$:
\begin{equation*}
T(N,\pi/4)=2\pi\sqrt{\frac{N}{g}}\Bigg[1+\sum_{n=1}^{\infty}\left(\frac{(2n)!}{2^{2\pi}(n!)^2}\right)^2\sin^{2n}\left(\frac{\pi/4}{2}\right)\Bigg].
\end{equation*} By running a Python simulation, we empirically verified this result and showed that, with the exception of a few outliers, the data reasonably fit the model for small $N$. When $N=1$ for a simple pendulum of unit length $1$ m and unit mass $1$ kg, the real time period is
\begin{equation*}
T(1,\pi/4)=2\pi\sqrt{\frac{1}{g}}\Bigg[1+\sum_{n=1}^{\infty}\left(\frac{(2n)!}{2^{2\pi}(n!)^2}\right)^2\sin^{2n}\left(\frac{\pi/4}{2}\right)\Bigg].
\end{equation*}

There are several improvements to the current model that would more closely represent a real physical system. Namely, we should relax the ideal constraints by considering extensible rods, non-uniform spacious mass bobs, friction due to air resistance, and a moving support. Likewise, systematic error can be minimized by limiting zero setting error and multiplier factor error in the measurement of the vertical equilibrium position. For initial oscillation amplitude $\theta_0=\pi/2$ rad, the system would be in unstable equilibrium. Furthermore, the system becomes exponentially unstable as $N$ increases. In fact, the implementation of the \code{Simulation()} method uses an iterative process that has exponential $O(m^N)$ run time for some $m\in\mathbb{N}$. However, this algorithm can be optimized to $\log$-run time $O(m\log N)$.

\newpage

\begin{appendices}
\section{Deriving the Euler--Lagrange Equations of Motion for an $N$-Body System}\label{sec:appendixa}
We derive the Euler--Lagrange equations of motion for the $N$-body planar system, namely $\frac{\partial L}{\partial\theta_i}=\frac{d}{dt}\left(\frac{\partial L}{\partial\dot\theta_i}\right)$. Recall the partial derivative of the Lagrangian with respect to the $\theta_i$ coordinate is:
\begin{equation}\label{leftlagrangian}
\begin{split}
\frac{\partial L}{\partial\theta_i}&=\frac{\partial}{\partial\theta_i}\sum_{j=1}^N\frac{1}{2}m_j(\dot x_j^2+\dot y_j^2)-gm_jy_j\\&
=\sum_{j=1}^N\frac{1}{2}m_j\left[2\dot x_j\frac{\partial \dot x_j}{\partial\theta_i}+2y_j\frac{\partial \dot y_j}{\partial\theta_i}\right]-gm_j\frac{\partial y_j}{\partial\theta_i}\\&=\sum_{j=1}^Nm_j\left[\dot x_j\frac{\partial \dot x_j}{\partial\theta_i}+\dot y_j\frac{\partial \dot y_j}{\partial\theta_i}-g\frac{\partial y_j}{\partial\theta_i}\right].
\end{split}
\end{equation}
Thus, we must compute the three partial derivatives $\partial\dot x_j/\partial\theta_i$, $\partial\dot y_j/\partial\theta_i$, and $\partial y_j/\partial\theta_i$. Then 
\begin{equation}
\begin{split}
\frac{\partial\dot x_j}{\partial\theta_i}&=\frac{\partial}{\partial\theta_i}\sum_{k=1}^jl_k\dot\theta_j\cos\theta_k\\&=
\begin{cases} 
      \frac{\partial}{\partial\theta_i}(l_1\dot\theta_1\cos\theta_1+\dots+l_i\dot\theta_i\cos\theta_i+\dots+l_j\dot\theta_j\cos\theta_j), & \text{if }i\le j, \\
      \frac{\partial}{\partial\theta_i}(l_1\dot\theta_1\cos\theta_1+\dots+l_j\dot\theta_j\cos\theta_j), & \text{if } i>j
   \end{cases}\\
 &=\begin{cases}
-l_i\dot\theta_i\sin\theta_i, & \text{if }i\le j\\
 0, & \text{if } i>j
 \end{cases}\\
&=-\xi_{ij}l_i\dot\theta_i\sin\theta_i
\end{split}
\end{equation} where
\begin{equation}
\xi_{ij}=\begin{cases}
1, & \text{if } i\le j\\
0, & \text{else.}
\end{cases}
\end{equation}
Similarly, the partial derivative of the $\dot y_j$ velocity with respect to the displacement angle $\theta_i$ is
\begin{equation}
\begin{split}
\frac{\partial\dot y_j}{\partial\theta_i}&=\frac{\partial}{\partial\theta_i}\sum_{k=1}^jl_k\dot\theta_k\sin\theta_k\\&=
\begin{cases} 
      \frac{\partial}{\partial\theta_i}(l_1\dot\theta_1\sin\theta_1+\dots+l_i\dot\theta_i\sin\theta_i+\dots+l_j\dot\theta_j\sin\theta_j), & \text{if }i\le j, \\
      \frac{\partial}{\partial\theta_i}(l_1\dot\theta_1\sin\theta_1+\dots+l_j\dot\theta_j\sin\theta_j), & \text{if } i>j
   \end{cases}\\
 &=\begin{cases}
l_i\dot\theta_i\cos\theta_i, & \text{if }i\le j\\
 0, & \text{if } i>j
 \end{cases}\\
&=\xi_{ij}l_i\dot\theta_i\cos\theta_i.
\end{split}
\end{equation} Lastly, the derivative of the vertical $y_j$ position with respect to the displacement angle $\theta_i$ is
\begin{equation}
\begin{split}
\frac{\partial y_j}{\partial\theta_i}&-=\frac{\partial}{\partial\theta_i}\sum_{k=1}^jl_k\cos\theta_k\\&=
\begin{cases} 
      \frac{\partial}{\partial\theta_i}(l_1\cos\theta_1+\dots+l_i\cos\theta_i+\dots+l_j\cos\theta_j), & \text{if }i\le j, \\
      \frac{\partial}{\partial\theta_i}(l_1\cos\theta_1+\dots+l_j\cos\theta_j), & \text{if } i>j
   \end{cases}\\
 &=-\begin{cases}
-l_i\sin\theta_i, & \text{if }i\le j\\
 0, & \text{if } i>j
 \end{cases}\\
&=\xi_{ij}l_i\sin\theta_i.
\end{split}
\end{equation}
It follows that Equation \ref{leftlagrangian} reduces to
\begin{equation}
\begin{split}
\frac{\partial L}{\partial\theta_i}&=\sum_{j=1}^Nm_j\left[\dot x_j\frac{\partial \dot x_j}{\partial\theta_i}+\dot y_j\frac{\partial \dot y_j}{\partial\theta_i}-g\frac{\partial y_j}{\partial\theta_i}\right]\\
&=\sum_{j=1}^Nm_j\left[\dot x_j(-\xi_{ij}l_i\dot\theta_i\sin\theta_i)+\dot y_j\xi_{ij}l_i\dot\theta_i\cos\theta_i-g\xi_{ij}l_i\sin\theta_i\right]\\
&=\sum_{j=1}^Nm_j\xi_{ij}l_i[(\dot y_j\dot\theta_i\cos\theta_i-\dot x_j\dot\theta_i\sin\theta_i)-g\sin\theta_i]\\
&=-\sum_{j=1}^N\xi_{ij}l_im_j\left[\dot\theta\sum_{k=1}^jl_j\dot\theta_k\sin(\theta_i-\theta_k)+g\sin\theta_i\right]\\
&=-l_i\left(\sum_{j=1}^{i-1}+\sum_{j=i}^N\right)m_j\left[\dot\theta_i\sum_{k=1}^jl_k\dot\theta_k\sin(\theta_i-\theta_k)+g\sin\theta_i\right]\\
&=-l_i\sum_{j=i}^Nm_j\left[\dot\theta_i\sum_{k=1}^jl_k\dot\theta_k\sin(\theta_i-\theta_k)+g\sin\theta_i\right]\\
&=-gl_i\sin\theta_i\sum_{j=i}^Nm_j-l_i\dot\theta_i\sum_{j=1}^Nl_j\dot\theta_j\sin(\theta_i-\theta_j)\sum_{k=\pi_{ij}}^Nm_k
\end{split}
\end{equation} where 
\begin{equation}
\pi_{ij}=\begin{cases}
i, & \text{if } j\le i\\
j, & \text{else.}
\end{cases}
\end{equation}

Likewise, we compute the partial derivative of the Lagrangian with respect to $i$-th angular velocity $\dot\theta_i$ to fully determine the Euler--Lagrange equations of motion:
\begin{equation}\label{rightlagrangian}
\begin{split}
\frac{\partial L}{\partial\dot\theta_i}=&\frac{1}{2}\sum_{j=1}^Nm_j\left[2\dot x_j\frac{\partial \dot x_j}{\partial\dot\theta_i}+2\dot y_j\frac{\partial\dot y_j}{\partial\dot\theta_i}\right]-g\sum_{i=1}^Nm_j\frac{\partial y_j}{\partial\dot\theta_i}\\
&=\sum_{i=1}^Nm_j\left[\dot x_j\frac{\partial\dot x_j}{\partial\dot\theta_i}+\dot y_j\frac{\partial\dot y_j}{\partial\dot\theta_i}-g\frac{\partial\dot y_j}{\partial\dot\theta_i}\right].
\end{split}
\end{equation} Recall that $\dot x_i=\sum_{k=1}^il_k\dot\theta_k\cos\theta_k$, $\dot y_i=\sum_{k=1}^il_k\dot\theta_k\sin\theta_k$, and $y_i=-\sum_{k=1}^il_k\cos\theta_k$. As such, 
\begin{equation}
\begin{split}
\frac{\partial\dot x_j}{\partial\dot\theta_i}&=\frac{\partial}{\partial\dot\theta_i}\sum_{k=1}^jl_k\dot\theta_k\cos\theta_k\\
&=\xi_{ij}l_i\cos\theta_i,\\
\frac{\partial \dot y_j}{\partial\dot\theta_i}&=\frac{\partial}{\partial\dot\theta_i}\sum_{k=1}^jl_k\dot\theta_k\sin\theta_k\\
&=\xi_{ij}l_i\sin\theta_i,
\end{split}
\end{equation} and
\begin{equation}
\frac{\partial y_j}{\partial\dot\theta_i}=-\frac{\partial}{\partial\theta_I}\sum_{k=1}^jl_k\cos\theta_k=0.
\end{equation} Thus, we substitute these partial derivatives in Equation \ref{rightlagrangian} to find
\begin{equation}
\begin{split}
\frac{\partial L}{\partial\dot\theta_i}&=\sum_{j=1}^N\xi_{ij}l_im_j[\dot y_j\sin\theta_i+\dot x_j\cos\theta_i]\\
&=\sum_{j=1}^N\xi_{ij}l_im_j\left[\sum_{k=1}^jl_j\dot\theta_k\sin\theta_k\sin\theta_i+\sum_{k=1}^jl_k\dot\theta_k\cos\theta_k\cos\theta_i\right]\\
&=\sum_{j=i}^Nl_im_j\left[\sum_{k=1}^jl_k\dot\theta_k(\sin\theta_k\sin\theta_i+\cos\theta_k\cos\theta_i)\right].
\end{split}
\end{equation} Since
\begin{equation*}
\begin{split}
\frac{\partial L}{\partial\dot\theta_i}&=\sum_{j=1}^Nm_jl_j\xi_{ij}[\dot y_j\sin\theta_i+\dot x_j\cos\theta_i]\\
&=l_i\sum_{j=i}^Nm_j(\dot x_j\cos\theta_i+\dot y_j\sin\theta_i),
\end{split}
\end{equation*} it follows that
\begin{equation}
\begin{split}
\frac{d}{dt}\left(\frac{\partial L}{\partial\dot\theta_i}\right)&=l_i\sum_{j=1}^Nm_j(\ddot x_j\cos\theta_i-\dot x_j\dot\theta_i\sin\theta_i+\ddot y_j\sin\theta_i+\dot y_j\dot\theta_i\cos\theta_i)\\
&=l_i\sum_{j=1}^Nm_j\Bigg[\sum_{k=1}^jl_k(\ddot\theta_k\cos\theta_k-\dot\theta_k^2\sin\theta_k)\cos\theta_i-\sum_{k=1}^jl_k\dot\theta_k\cos\theta_k\dot\theta_i\sin\theta_i\\
&+\sum_{k=1}^jl_k(\ddot\theta_k\sin\theta_k+\dot\theta_k^2\cos\theta_k)\sin\theta_i+\sum_{k=1}^kl_k\dot\theta_k\sin\theta_k\dot\theta_i\cos\theta_i\Bigg]\\
&=l_i\sum_{j=i}^Nm_j\left[\sum_{k=1}^jl_k(\ddot\theta_j\cos(\theta_i-\theta_j)+(\dot\theta_k^2-\dot\theta_i\dot\theta_k)\sin(\theta_i-\theta_k))\right]\\
&=l_i\sum_{j=1}^Nl_j\left[\ddot\theta_j\cos(\theta_i-\theta_j)+(\dot\theta_j^2-\dot\theta_i\dot\theta_j)\sin(\theta_i-\theta_j)\right]\sum_{k=\pi_{ij}}^Nm_k.
\end{split}
\end{equation} 

\section{Linearizing the Non-Linear Euler--Lagrange Equations of Motion}\label{sec:appendixlin}
The Maclaurin series developments of sine and cosine 
\begin{equation}
\begin{split}
&\cos\theta =\sum_{n=0}^{\infty}(-1)^n\frac{\theta^{2n}}{2n!}\approx \theta,\\
&\sin\theta=\sum_{n=0}^{\infty}(-1)^n\frac{\theta^{2n+1}}{(2n+1)!} \approx 1
\end{split}
\end{equation}
for small $\theta$ imply that $\sin\theta_i\approx \theta_i$, $\sin(\theta_i-\theta_j)\approx \theta_i-\theta_j,$ and $\cos(\theta_i-\theta_j)=1-\frac{(\theta_i-\theta_j)^2}{2}\approx 1$. Thus, the perturbative expansions dictate, in tandem, that the equations of motion corresponding to the linearized Lagrangian are 
\begin{equation}\label{lineom}
\sum_{j=1}^Nl_j\sum_{k=\pi_{ij}}^Nm_k\ddot\theta_j=-g\theta_i\sum_{j=i}^Nm_j-\sum_{j=1}^Nl_j\dot\theta_j^2(\theta_i-\theta_j)\sum_{l=\pi_{ij}}^Nm_k.
\end{equation} 
We may re-write this as
\begin{equation}
\sum_{j=1}^NM_{ij}\ddot\theta_j=\sum_{j=1}^NK_{ij}
\end{equation} so by row: $M_{ij}\ddot\theta_j=K_{ij}$ where $M_{ij}=l_j\sum_{k=\pi_{ij}}^Nm_k$ and, for the approximation $\sum_{j=1}^Nm_j\approx \sum_{j=i}^Nm_j$, we have $K_{ij}\approx -g\theta_im_j-l_j\dot\theta_j^2(\theta_i-\theta_j)\sum_{k=\pi_{ij}}^Nm_k\approx l_j\dot\theta_j^2\theta_j\sum_{k=\pi_{ij}}^Nm_k+g\theta_jm_j=-\theta_j\left[l_j\dot\theta_j^2\sum_{k=\pi_{ij}}^Nm_k+gm_j\right]$. Note, we linearize $\theta_i\approx\theta_j$ to eliminate non-linear behavior in the dynamical system.
In matrix notation, we may write the $(i,j)$ entry as:
\begin{equation}
\left[l_j\sum_{k=\pi_{ij}}^Nm_k\right]\ddot\theta_j+\left[l_j\sum_{k=\pi_{ij}}^Nm_k\right]\dot\theta_j^2+\left[gm_j\right]\theta_j=\mathbf{0}.
\end{equation}

\section{Fourth Order Runge-Kutta Method}\label{sec:appendixrk4}
The system in Equation \ref{system} may be re-written by observing 
\begin{equation}
\begin{bmatrix}
           \dot{\bm{\alpha}}_1(t)\\
           \dot{\bm{\alpha}}_2(t)
         \end{bmatrix}:=
         \begin{bmatrix}
          \mathbf{F}_1(t,\bm{\alpha}_2)\\
           \mathbf{F}_2(t,\bm{\alpha}_1,\bm{\alpha}_2)
         \end{bmatrix}.   
\end{equation} Define the matrices $\bm{\alpha}=(\bm{\alpha}_1,\bm{\alpha}_2)^T$ such that:
\begin{equation}
\frac{d\bm{\alpha}}{dt}=\mathbf{F}(t,\bm{\alpha}).
\end{equation} This is an initial value problem with $\bm{\alpha}(t_0)=\bm{\alpha}_0.$ We find $\bm{\alpha}_{n+1}$, the RK4 approximation of $\bm{\alpha}(t_{n+1})$, which is uniquely determined by the previous value $\bm{\alpha}_{n}$ plus a weighted average of four vectors which are given by the product of the step size $\tau$ and the estimated value of the slope field in phase space given by $\mathbf{F}(t,\bm{\alpha})$. Then RK4 is implemented by defining the following matrices \cite{math24}:
\begin{equation}
\begin{split}
&\mathbf{K}_1=\tau\mathbf{F}(t_n,\bm{\alpha}(t_n)),\\
&\mathbf{K}_2=\tau\mathbf{F}\left(t_n+\frac{1}{2}\tau,\bm{\alpha}(t_n)+\frac{1}{2}\mathbf{K}_1\right)\\
&\mathbf{K}_3=\tau\mathbf{F}\left(t_n+\frac{1}{2}\tau,\bm{\alpha}(t_n)+\frac{1}{2}\mathbf{K}_2\right)\\
&\mathbf{K}_4=\tau\mathbf{F}\left(t_n+\tau,\bm{\alpha}(t_n)+\mathbf{K}_3\right)
\end{split}
\end{equation} 
whereby the next vector in the node (in the neighborhood) of the integral curve in phase space is:
\begin{equation}
\bm{\alpha}(t_{n+1})=\bm{\alpha}(t_n)+\frac{1}{6}(\mathbf{K}_1+2\mathbf{K}_2+2\mathbf{K}_3+\mathbf{Y}_4)
\end{equation} with two initial conditions specified by $\bm{\alpha}(t_0).$
\section{Kane's Method Python Implementation for Euler--Lagrange Equations}\label{sec:appendix}

The following Kane's method integrator for equations of motion is due to VanderPlas \cite{jakevdp} and Gede et al. \cite{gede}.
\fontsize{7}{9}\selectfont
\begin{lstlisting}[language=Python, caption=Kane's Method Equations of Motion Integrator]
import matplotlib.pyplot as plt
import numpy as np
import pandas as pd

from sympy import symbols
from sympy.physics import mechanics

from sympy import Dummy, lambdify
from scipy.integrate import odeint


def integrate_pendulum(n, times,
                       initial_positions,
                       initial_velocities=0,
                       lengths=None, masses=1):
    """Integrate a multi-pendulum with `n` sections"""
    #-------------------------------------------------
    # Step 1: construct the pendulum model
    
    # Generalized coordinates and velocities
    # (in this case, angular positions & velocities of each mass) 
    q = mechanics.dynamicsymbols('q:{0}'.format(n))
    u = mechanics.dynamicsymbols('u:{0}'.format(n))

    # mass and length
    m = symbols('m:{0}'.format(n))
    l = symbols('l:{0}'.format(n))

    # gravity and time symbols
    g, t = symbols('g,t')
    
    #--------------------------------------------------
    # Step 2: build the model using Kane's Method

    # Create pivot point reference frame
    A = mechanics.ReferenceFrame('A')
    P = mechanics.Point('P')
    P.set_vel(A, 0)

    # lists to hold particles, forces, and kinetic ODEs
    # for each pendulum in the chain
    particles = []
    forces = []
    kinetic_odes = []

    for i in range(n):
        # Create a reference frame following the i^th mass
        Ai = A.orientnew('A' + str(i), 'Axis', [q[i], A.z])
        Ai.set_ang_vel(A, u[i] * A.z)

        # Create a point in this reference frame
        Pi = P.locatenew('P' + str(i), l[i] * Ai.x)
        Pi.v2pt_theory(P, A, Ai)

        # Create a new particle of mass m[i] at this point
        Pai = mechanics.Particle('Pa' + str(i), Pi, m[i])
        particles.append(Pai)

        # Set forces & compute kinematic ODE
        forces.append((Pi, m[i] * g * A.x))
        kinetic_odes.append(q[i].diff(t) - u[i])

        P = Pi

    # Generate equations of motion
    KM = mechanics.KanesMethod(A, q_ind=q, u_ind=u,
                               kd_eqs=kinetic_odes)
    fr, fr_star = KM.kanes_equations(particles, forces)
    
    #-----------------------------------------------------
    # Step 3: numerically evaluate equations and integrate

    # initial positions and velocities ? assumed to be given in degrees
    y0 = np.deg2rad(np.concatenate([np.broadcast_to(initial_positions, n),
                                    np.broadcast_to(initial_velocities, n)]))
        
    # lengths and masses
    if lengths is None:
        lengths = np.ones(n) / n
    lengths = np.broadcast_to(lengths, n)
    masses = np.broadcast_to(masses, n)

    # Fixed parameters: gravitational constant, lengths, and masses
    parameters = [g] + list(l) + list(m)
    parameter_vals = [9.81] + list(lengths) + list(masses)

    # define symbols for unknown parameters
    unknowns = [Dummy() for i in q + u]
    unknown_dict = dict(zip(q + u, unknowns))
    kds = KM.kindiffdict()

    # substitute unknown symbols for qdot terms
    mm_sym = KM.mass_matrix_full.subs(kds).subs(unknown_dict)
    fo_sym = KM.forcing_full.subs(kds).subs(unknown_dict)

    # create functions for numerical calculation 
    mm_func = lambdify(unknowns + parameters, mm_sym)
    fo_func = lambdify(unknowns + parameters, fo_sym)

    # function which computes the derivatives of parameters
    def gradient(y, t, args):
        vals = np.concatenate((y, args))
        sol = np.linalg.solve(mm_func(*vals), fo_func(*vals))
        return np.array(sol).T[0]

    # ODE integration
    return odeint(gradient, y0, times, args=(parameter_vals,)) 
\end{lstlisting}

We extract the Cartesian $(x,y)$ coordinates of each pendulum by implementing the following method:

\begin{lstlisting}[language=Python, caption=Extracting Cartesian coordinates]
def get_xy_coords(p, lengths=None):
    """Get (x, y) coordinates from generalized coordinates p"""
    p = np.atleast_2d(p)
    n = p.shape[1] // 2
    if lengths is None:
        lengths = np.ones(n) / n
    zeros = np.zeros(p.shape[0])[:, None]
    x = np.hstack([zeros, lengths * np.sin(p[:, :n])])
    y = np.hstack([zeros, -lengths * np.cos(p[:, :n])])
    return np.cumsum(x, 1), np.cumsum(y, 1)
\end{lstlisting}
\fontsize{10}{12}\selectfont

Thus, we implement a \code{Simulation()} method which will apply the aforementioned procedure to determine the average oscillation period for each pendulum bob and then take an average over all masses. When we call \code{Simulation()}, the local method \code{system_period_n(i)} returns the average time period of all $i$ pendulums after calling the \code{get_timeperiod_y(j)} method to get the average time period of the $j$-th pendulum bob for all times. We call this method for $i=5$, $10$, $20$, $100$. The example below has static field number of pendulums \code{numberofps=20}, but we amend this depending on the number $N$ we would like. Hence, we generate the \code{sequence_of_periods} for each trial after calling the \code{Simulation()} method which perturbs the initial angular positions \code{initial_positions=45}, measured in degrees, to a random float in the range $[0,1]$ by importing the \code{random} library. The $1\degree$ random addition accounts for the uncertainty error of measurement for a degree.
\fontsize{7}{9}\selectfont
\begin{lstlisting}[language=Python, caption=Simulated Time Period Method]
import math
import numpy as np
import matplotlib.pyplot as plt
import random

def Simulation():
    numberofps = 20
    n = numberofps + 1
    timeintercepts = []
    nperiod = []
    sequence_of_periods = []
    deltheta_0 = random.randint(0,1)
    for i in range(1,n + 1):
        t = np.linspace(0, 10, 1000)
        p = integrate_pendulum(i, times=t, initial_positions=45 + deltheta_0)
        x, y = get_xy_coords(p)
        #plt.plot(x, y);
        #x, y has first column of zeros
        r,s = np.shape(y)
        def get_r():
            return r
        def theta(j):
            theta_j = []
            for i in range(0, r):
                theta_j.append(math.acos(abs(y[i][j-1]-y[i][j]/1)))                					    theta_j.append(math.acos(abs(y[i][j-1]-y[i][j]/1)))
            timenew = [i for i in range(1,r + 1)]
            graph_j = pd.Series(data=theta_j, index=timenew)
            return pd.Series(data=np.gradient(graph_j.values), index=graph_j.index) #returns omega_j
        def computeperiod():
            series = []
            for j in range(1,s):
                series.append(2 * math.pi/(theta(j).mean()))
               ## return 2 * math.pi/(theta(j).mean())
            numberline = [i for i in range(1,s)] #Here s=n+1
            timeperiod = pd.Series(data=series, index=numberline)
            return abs(timeperiod.mean())
        def get_timeperiod_y(j):
            #get y_j
            y_j = []
            for i in range(0, r):
                y_j.append(y[i][j])
            z_j = np.gradient(y_j)
            times = []
            #CHECKING FOR CRITICAL POINTS
            for k in range(len(z_j)):
                if np.sign(z_j[k]) != np.sign(z_j[k-1]) and k-1 not in times:
                    times.append(k*1/100)
            #INSERT CONDITION CHECKING FOR LOOP TRAJECTORIES, I.E. y(t_i) = y(t_j)
            avg_j = []
            tmp = []
            tmp1 = []
            min_1 = 0
            min_2 = 0
            for p in range(len(y_j)):
                for q in range(len(y_j)):
                    if(abs(y_j[p]-y_j[q]) < 0.01 and p < q):
                        tmp1.append(p)
                        tmp1.append(q)
                        min_1 = min(tmp1)
                        tmp1.remove(min_1)
                        min_2 = min(tmp1)
                        tmp.append(abs(min_2-min_1)*1/100) 
            avg_j.append(abs(np.mean(tmp)))
            for i in range(0,len(times)-2):
                avg_j.append(times[i+2]-times[i])
            return abs(np.mean(avg_j))
        def system_period_n(n):
            allperiods = []
            for j in range(1, n + 1):
                allperiods.append(1/(get_timeperiod_y(j)))
            periodresult_n = np.mean(allperiods)
            return periodresult_n
        sequence_of_periods.append(system_period_n(i))
    return sequence_of_periods

sequence_of_periods_1 = Simulation()
sequence_of_periods_2 = Simulation()
sequence_of_periods_3 = Simulation()
\end{lstlisting}

\fontsize{10}{12}\selectfont
Finally, we graph the pseudo-period as a function of the number of pendulums $N$, which we write in the code as \code{n}, for each trial along with error bars due to higher-order corrections of Equation \ref{periodexact}.

\fontsize{7}{9}\selectfont
\begin{lstlisting}[language=Python, caption=Visualizing the pseudo-time period]

x_1 = [sequence_of_periods_0[i] for i in range(0,n+1)]
x_2 = [sequence_of_periods_2[i] for i in range(0,n+1)]
x_3 = [sequence_of_periods_3[i] for i in range(0,n+!)]

def period(n):
    if n != 0:
        m = [1] * n
        l = [1] * n
        M = sum(m)
        num = 2 * math.pi * n
        dem = 0
        g = 9.8
        for i in range(n):
            dem = dem + math.sqrt(g * m[i]/(l[i]*M))
        return num/dem 
    elif n == 0:
        return 0
x = np.arange(0, n + 1,1)
y = np.array([period(i) for i in x])
yerr_val = (0.0025 + 2/135+4/135+6/135+8/135+10/135+12/135+14/135)*2*(math.pi)*np.sqrt(x/9.8)
plt.errorbar(x, y, yerr=yerr_val)

# Data
df = pd.DataFrame({'number': range(0,n + 1), 'Model': y, 'Trial 1': x_1, 'Trial 2': x_2, 'Trial 3': x_3})
 
# multiple line plot
plt.plot( 'number', 'Model', data=df, marker='', color='skyblue', linewidth=2)
plt.plot( 'number', 'Trial 1', data=df, marker='', color='black', linewidth=2, linestyle = 'dashed')
plt.plot( 'number', 'Trial 2', data=df, marker='', color='green', linewidth=2, linestyle = 'dashed')
plt.plot( 'number', 'Trial 3', data=df, marker='', color='red', linewidth=2, linestyle = 'dashed')
plt.legend()
plt.title("Pseudo-time period for system of n-pendulums")
plt.xlabel('Number of pendulums (n)')
plt.ylabel('Pseudo-period T(n)/s')
\end{lstlisting} 

\end{appendices}
\fontsize{10}{12}\selectfont


\newpage

\newpage
\begin{bibdiv}
\begin{biblist}

\bib{shankar}{book}{
  title={Fundamentals of physics: mechanics, relativity, and thermodynamics},
  author={Shankar, R.},
  date={2014},
  publisher={Yale University Press}
  address = {New Haven, Connecticut, USA}
}

\bib{nelson}{article}{
title = {The pendulum--rich physics from a simple system},
author = {Nelson, A.~R.}
journal = {Am. J. Phys},
   volume = {54},
     date = {1986},
    pages = {112-121},
}

\bib{sympy}{article}{
title = {Kane's method in physics/mechanics},
author = {SymPy Development Team}, 
year = {2018},
url = {https://docs.sympy.org/latest/modules/physics/mechanics/kane.html},
eprint = {https://docs.sympy.org/latest/modules/physics/mechanics/kane.html}
}

\bib{math24}{article}{
title = {Double pendulum},
author = {Math24.net}, 
year = {2018},
url = {https://www.math24.net/double-pendulum/},
eprint = {https://www.math24.net/double-pendulum/}
}

\bib{jakevdp}{article}{
title = {Double pendulum},
author = {VanderPlas, J.}, 
year = {2017},
url = {https://jakevdp.github.io/blog/2017/03/08/triple-pendulum-chaos/},
eprint = {https://jakevdp.github.io/blog/2017/03/08/triple-pendulum-chaos/}
}

\bib{gede}{article}{
title = {Constrained multibody dynamics with Python: from symbolic equation generation to publication},
author = {Gede, G.}, 
author = {Peterson, D~L.}
author = {Nanjangud, A~S.}
author = {Moore, J~K.}
author = {Hubbard, M.}
journal = {Proceedings of the ASME 2013 IDETC/CIE}
year = {2013},
}

\end{biblist}
\end{bibdiv}
\end{document}